\documentclass[twoside]{irmaems}
\usepackage{amssymb} 
\usepackage{amsmath}
\usepackage{latexsym}

\newcommand{\nb}[1]{#1\nobreakdash-}

\newcommand\ds\displaystyle

\theoremstyle{plain}
\newtheorem{theorem}{Theorem}[section]

\newtheorem{lemma}[theorem]{Lemma}

\DeclareMathOperator{\Fix}{Fix}
\DeclareMathOperator{\SL}{SL}
\DeclareMathOperator{\PSL}{PSL} 
\DeclareMathOperator{\Homeo}{Homeo}

\DeclareMathOperator{\Isom}{Isom}

\DeclareMathOperator{\Stab}{Stab}

\DeclareMathOperator\kernel{kernel}

\DeclareMathOperator\Hull{Hull}

\newcommand\R{{\mathbf R}}
\newcommand\reals{\R}

\newcommand\hyp{\mathbf{H}}
\newcommand\Complex{{\mathbf C}}
\newcommand\complex{\Complex}
\newcommand\Z{{\mathbf Z}}
\newcommand\integers{\Z}

\renewcommand\H{{\mathcal{H}}}

\newcommand\infinity{\infty}
\newcommand\bndry{\partial}
\newcommand{\bdy}{\bndry}
\newcommand{\from}{\colon}

\newcommand\suchthat{\bigm|}
\newcommand\inverse{{-1}}
\newcommand\inv{\inverse}

\newcommand\union{\cup}

\newcommand\absvalue[1]{\left| #1 \right|}
\newcommand\abs[1]{\absvalue{#1}}

\newcommand\wt\widetilde

\newcommand\intersect{\cap}

\DeclareMathOperator\QD{QD}

\newcommand\restrict{\bigm|}
\newcommand\subgroup{<}

\newcommand\Teichmuller{Teichm\"uller}
\newcommand\Poincare{Poincar\'e}

\newcommand\MCG{{\mathcal M \mathcal C \mathcal G}}

\newcommand\cross{\times}

\newcommand\M{{\cal M}}
\newcommand\Mod\M

\newcommand\PMF{\mathcal{PMF}}
\newcommand\MF{{\cal MF}}
\newcommand\GL{{\cal GL}}

\newcommand\Union\bigcup
\newcommand\C{\mathcal C}

\newcommand\Teich{\T}
\newcommand\T{{\cal T}}
\newcommand\<\langle
\renewcommand\>\rangle

\newcommand\carries\succ
\newcommand\carriedby\prec
\newcommand\splitsto\carries
\newcommand\collapsesto\carriedby

\newcommand\fol{{\mathcal F}}

\newcommand\geodesic[2]{\overleftrightarrow{(#1,#2)}}

\newcommand\CC{\mathcal{CC}}
\newcommand\ceiling[1]{\lceil #1 \rceil}

\newcommand\Papa{Papadopoulos}
\newcommand\McPapa{McCarthy--\Papa}

\begin{document}

\title{Geometric survey of subgroups \\ of mapping class groups}

\author{Lee Mosher\thanks{
To appear in Handbook of Teichm\"uller Theory, Volume 1, ed. A. Papadopolous, European Math. Soc. 2007.
Work partially supported by NSF Grant No.~0405979
}}

\address{
Dept. of Mathematics and Computer Science\\
Rutgers University\\
Newark, NJ 07102\\
email:\,\tt{mosher@andromeda.rutgers.edu}
}

\maketitle

\begin{abstract}
The theme of this survey is that subgroups of the mapping class group $\MCG$ of a finite type surface $S$ can be studied via the geometric/dynamical properties of their action on the Thurston compactification of \Teichmuller\ space $\overline\Teich= \Teich \union \PMF$, just as discrete subgroups of the isometries of hyperbolic space $\hyp^n$ can be studied via their action on compactified hyperbolic space $\overline\hyp^n = \hyp^n \union S^{n-1}_\infinity$.
\end{abstract}



\tableofcontents

\section{Introduction}
Let $S$ be a finite type surface, obtained from a closed oriented surface of genus $g \ge 0$ by removing $p \ge 0$ punctures. The (orientation preserving) \emph{mapping class group} of $S$\index{mapping class group} is denoted $\MCG=\MCG(S) = \Homeo_+(S) / \Homeo_0(S)$ where $\Homeo_+(S)$ is the group of orientation preserving homeomoprhisms and $\Homeo_0(S)$ is the normal subgroup of homeomorphisms isotopic to the identity. The \emph{\Teichmuller\ space} of $S$, denoted $\Teich=\Teich(S)$, is the space of conformal structures on $S$ modulo isotopy, or equivalently the space of complete, finite area hyperbolic structures on $S$ modulo isotopy; the phrase ``modulo isotopy'' means ``modulo the action of $\Homeo_0(S)$''. Thurston's boundary of \Teichmuller\ space, denoted $\PMF=\PMF(S)$, is the space of projective Whitehead classes of measured foliations on $S$. There is a natural topology on the set $\overline\Teich = \Teich \union \PMF$ which makes it homeomorphic to a ball of dimension $6g-g+2p$ with interior $\Teich$ and boundary $\PMF$. The group $\MCG$ acts naturally on $\Teich$ and on $\PMF$, extending to an action by homeomorphisms on~$\overline\Teich$.

The theme of this survey is that subgroups of $\MCG$ can be studied via the geometric/dynamical properties of their action on $\overline\Teich= \Teich \union \PMF$, just as discrete subgroups of the isometries of hyperbolic space $\hyp^n$ can be studied via their action on $\overline\hyp^n = \hyp^n \union S^{n-1}_\infinity$. Recent problem lists with this theme in mind include \cite{Mosher:ExtensionProblems} and \cite{Reid:SurfaceSubgroups}. We will survey several topics exemplifying this theme, in the following sections:
\begin{description}
\item[\S \ref{SectionClassification}] The Thurston--Bers classification of mapping classes of $\MCG$.
\item[\S \ref{SectionTits}] The Tits alternative.
\item[\S \ref{SectionLimits}] Limit sets and domains of discontinuity.
\item[\S \ref{SectionConvexCocompact}] Convex cocompact subgroups.
\item[\S \ref{SectionStabilizers}] Stabilizers of \Teichmuller\ discs.
\item[\S \ref{SectionCombinations}] Combination theorems and surface subgroups.
\end{description}
Here is a brief description of each of these topics, to be filled out in the rest of the paper.

Thurston \cite{Thurston:surfaces} and Bers \cite{Bers:ThurstonTheorem} classified individual elements of $\MCG$, or equivalently cyclic subgroups of $\MCG$, via the geometry and dynamics of their actions on $\overline\Teich$. Finite order mapping classes have fixed points in $\Teich$. 
Infinite order irreducible mapping classes have positive translation distance in $\T$, which is minimized along a unique \Teichmuller\ geodesic, and they act on $\overline\T$ with source-sink dynamics. Finite order mapping classes have fixed points in $\T$, at which the translation distance of zero is minimized.  And infinite order reducible mapping classes have translation distance in $\T$ whose infimum is not realized.

The \emph{Tits alternative}, originally formulated and proved for lattices in semisimple Lie groups, was also formulated and proved for subgroups of $\MCG$ independently by McCarthy \cite{McCarthy:Tits} and Ivanov \cite{Ivanov:subgroups}. The Tits alternative says that each subgroup of $\MCG$ has either an abelian subgroup of finite index, or a free subgroup of rank $\ge 2$.

\emph{Limit sets and domains of discontinuity} of subgroups of $\MCG$ were studied by \McPapa\ \cite{McCarthyPapa:dynamics}, in analogy to limit sets and domains of discontinuity for discrete subgroups of $\Isom(\hyp^n)$.

\emph{Convex cocompact subgroups} of $\MCG$, including Schottky subgroups of $\MCG$, were originally studied by Farb and Mosher \cite{FarbMosher:quasiconvex} and further developed by Kent--Leininger \cite{KentLeininger:Shadows} and Hamenstadt \cite{Hamenstadt:extensions}, in analogy to convex cocompact subgroups and Schottky subgroups of $\Isom(\hyp^n)$.

\emph{\Teichmuller\ discs} in $\Teich$ are isometrically embedded copies of the hyperbolic plane, and their stabilizers are analogues in $\MCG$ of Fuchsian subgroups of $\Isom(\hyp^n)$, which are stabilizers of isometrically embedded copies of $\hyp^2$ in $\hyp^n$.

The \emph{Leininger--Reid combination theorem} \cite{LeiningerReid:combination} gives a method for building closed surface subgroups of $\MCG$, by combining certain \Teichmuller\ disc stabilizer subgroups, in analogy to the Maskit combination theorem and related results for building discrete subgroups of $\Isom(\hyp^n)$ by combining simpler discrete subgroups. 

Because of our focus on the geometric/dynamic properties of subgroups of $\MCG$, we have nothing to say about the very rich and interesting algebraic tools that are used to study subgroups, particularly the Torelli subgroup and the lower central series of $\MCG$.

\section{The action of $\MCG$ on compactified \Teichmuller\ space}

We start by establishing what will be the basic notions and notations for this survey. References for this material include \cite{ImayoshiTaniguchi} and \cite{FLP}.

\paragraph{\Teichmuller\ space and the Thurston boundary.}
An \emph{essential curve}\index{essential curve} on $S$ is a simple closed curve which bounds neither a disc nor a once-punctured disc. The isotopy class of an essential curve $c$ is denoted $[c]$. Let $\C$ be the set of isotopy classes of essential curves on $S$. An \emph{essential curve family} is a finite set of pairwise disjoint essential curves, no two of which are isotopic. The \emph{\Teichmuller\ space} of $S$, denoted $\Teich(S)$, is the set of complete, finite area hyperbolic structures on $S$, or equivalently the set of conformal structures with removable singularities at the punctures, modulo isotopy. Isotopy classes, in $\C$ or in $\Teich$ or in other contexts, will be denoted formally with square brackets $[\cdot]$, although informally these brackets will often be dropped, as in the next sentence. There is an embedding $\Teich(S) \mapsto \reals^\C$ defined by assigning to each $\sigma \in \Teich$ and $c \in \C$ the minimal length of curves isotopic to $c$ in the hyperbolic structure $\sigma$. The topology on $\Teich(S)$ is defined so that the map $\Teich(S) \mapsto \reals^\C$ is a homeomorphism onto its image. The group $\MCG(S)$ acts naturally on $\C$ and on $\Teich(S)$, and so the embedding $\Teich(S) \mapsto \reals^\C$ is $\MCG$-equivariant. 

A \emph{quadratic differential} on $\sigma \in \Teich$ is, formally, a holomorphic section of the symmetric square of the cotangent bundle of $\sigma$, whose integral is finite. Such sections can be added and multiplied by complex scalars, and so the quadratic differentials on $\sigma$ form a vector space $\QD_\sigma$ over the complex numbers. Informally, to each conformal coordinate $z$ for $\sigma$, a quadratic differential $\theta$ assigns an expression $f(z) dz^2$, where $f(z)$ is holomorphic and such that the expression behaves correctly under coordinate change. For $p \in S$ outside of a finite singular set, $\theta$ has a local coordinate $z$ in which its expression is $dz^2$ and $z(p)=0$, and $z$ is unique up to multiplication by $\pm 1$. On each point of the singular set, including each puncture, $\theta$ has a local coordinate $z$ in which its expression is $z^{k-2} dz^2$ for some integer $k$, where $k \ge 3$ if $p$ is not a puncture and $k \ge 1$ if $p$ is a puncture, and this coordinate $z$ is well-defined up to multiplication by a $k^{\text{th}}$ root of unity. The \emph{horizontal and vertical measured foliations} of $\theta$ are defined in each regular canonical coordinate $z=x+iy$ as follows: the horizontal foliation has leaves on which $y$ is constant, with transverse measure $\abs{dx}$; the vertical foliation has leaves on which $x$ is constant, with transverse measure $\abs{dy}$. At a singularity $p$ near which $\theta$ has the canonical expression $z^{k-2} dz^2$ we say that the horizontal and vertical measured foliations have \emph{$k$-pronged singularities} at $p$.

\Teichmuller's theorem says that for any $\sigma,\sigma' \in \Teich$ there exist unique choice of quadratic differentials $\theta$ on $\sigma$ called the \emph{initial quadratic differential}, $\theta'$ on $\sigma'$ called the \emph{terminal quadratic differential}, and a number $d \ge 0$, such that if $z=x+iy$ is a regular local coordinate for $\theta$ then (up to isotopy of $\sigma'$) the expression $z' = e^d x + i e^{-d} y$ defines a regular local coordinate for $\theta'$. This number $d$ is called the \emph{\Teichmuller\ distance} between $\sigma$ and $\sigma'$, and it defines a metric on $\Teich(S)$ which gives the same topology as that defined above.

Composing the embedding $\Teich(S) \mapsto \reals^\C$ with the projectivization $\reals^\C \mapsto P\reals^\C$, the composition $\Teich(S) \mapsto P\reals^\C$ is still an embedding, and its image is precompact. The closure of the image minus the image defines \emph{Thurston's boundary} for $\Teich(S)$, and this boundary is characterized as follows. A \emph{measured foliation}\index{measured foliation} on $S$ is a foliation with a finite singular set, equipped with a positive transverse Borel measure, so that at each singularity $p$ including punctures there exists $k \ge 1$ such that the singularity is locally modelled on a $k$-pronged singularity of $z^{k-2} dz^2$, where $k \ge 3$ if $p$ is not a puncture. There is an equivalence relation on measured foliations generated by isotopy and by the \emph{Whitehead move}, a move which collapses to a point each leaf segment whose endpoints are a pair of singularities at least one of which is not a puncture. The equivalence classes are called \emph{Whitehead classes}, and the set of Whitehead classes of measured foliations is denoted $\MF$ (again we use $[\cdot]$ to denote Whitehead class, except that the notation is often dropped). There is an embedding $\MF \mapsto \reals^\C$ defined by assigning to each $\fol \in \MF$ and $c \in \C$ the number $\<\fol,c\> \ge 0$ which is the minimal integral of curves homotopic to $c$ with respect to the transverse measure on $\fol$. The topology on $\MF$ is chosen so that the map $\MF \mapsto \reals^\C$ is a homeomorphism onto its image. When the transverse measure on $\fol$ is multiplied by a number $r>0$ we denote the result $r\fol$, and this generates the relation of \emph{projective equivalence} on $\MF$. The space of projective equivalence classes is denoted $\PMF$, and the resulting map $\PMF \mapsto P\reals^\C$ is an embedding. 

Thurston's compactification theorem says that the disjoint union of $\Teich$ and $\PMF$ embeds into $P\reals^\C$, the image is a closed ball of dimension $6g-6+2p$, the interior of this ball is the image of $\Teich$, and its boundary sphere is the image of $\PMF$. There is a natural action of $\MCG$ on $\PMF$ with respect to which the embedding $\PMF \to P\reals^\C$ is equivariant, and hence $\MCG$ acts homeomorphically on the ball $\overline\Teich = \Teich \union \PMF$. 

Consider an essential closed curve $c \in \C$. The \emph{enlargement} of $c$ is a measured foliation $\fol_c$ well-defined up to Whitehead equivalence as follows. First one chooses a spine for $S-c$ containing all punctures. The complement of this spine is an open annulus with core $c$. Then one foliates this annulus by curves isotopic to $c$, which together with the spine defines $\fol_c$ as a singular foliation. The transverse measure on $\fol_c$ is defined so that a curve in the annulus from boundary to boundary transverse to $\fol_c$ has total transverse measure~1. This construction generalizes to the \emph{enlargement of a weighted multicurve}, a formal sum $a_1 c_1 + \cdots + a_k c_k$ where $\{c_1,\ldots,c_k\}$ is an essential curve system and $a_1,\ldots,a_k > 0$; in this case one chooses a spine for $S-(c_1 \union \cdots \union c_k)$ containing all punctures, whose complement is a union of open annuli with cores $c_1,\ldots,c_k$ respectively; then one foliates these annuli by closed curves parallel to the cores; then one assigns total transverse measures $a_1,\ldots,a_k$, respectively, to these foliated annuli. The set of points in $\PMF$ represented by enlargements of weighted multicurves is dense. In fact, the set of points in $\MF$ represented by enlargements of essential closed curves is dense.

The intersection pairing $\MF \cross \C \to [0,+\infinity)$ taking $\fol,c$ to $\<\fol,c\>$ extends to a continuous intersection pairing $\MF \cross \MF \to [0,+\infinity)$, called the \emph{intersection number} of $\fol,\fol'$; uniqueness of the extension follows from denseness of $\C$ in $\MF$. 

A pair of measured foliations $\fol,\fol'$ is said to \emph{jointly fill} the surface $S$ if for every $c \in \C$ either $\<\fol,c\>$ or $\<\fol',c\>$ is nonzero, and $\<\fol,\fol'\>$ is also nonzero. This happens if and only if there exists a conformal structure $\sigma \in \Teich$ and a quadratic differential $\theta$ on $\sigma$ such that $\fol,\fol'$ are Whitehead equivalent to the horizontal and vertical measured foliations of $\theta$; moreover, $\sigma$ and $\theta$ are uniquely determined by $\fol,\fol'$.

The set of geodesic lines in $\Teich$ can be parameterized as follows. Choose measured foliations $\fol,\fol'$ which jointly fill. For each $r \in \reals$, the measured foliations $e^r \fol, e^{-r}\fol'$ jointly fill and determined a conformal structure $\sigma_r$. The map $r \mapsto \sigma_r$ is a geodesic embedding $\reals \to \Teich$. The image of this geodesic depends only on the projective classes $P[\fol]$, $P[\fol']$, and we say that this is the \emph{\Teichmuller\ geodesic} with directions $P[\fol]$, $P[\fol']$. 

The set of geodesic rays in the \Teichmuller\ metric can be parameterized as follows. First we use a theorem of Hubbard and Masur \cite{HubbardMasur:qd} which says that for each $\sigma \in \Teich$ and $P[\fol] \in \PMF$ there is up to real scalar multiple a unique $\theta \in \QD_\sigma$ whose horizontal measured foliation is in the given projective class $P[\fol]$. For each $d>0$ we can then transform each regular canonical coordinate $z=x+iy$ of $\theta$ to get a coordinate $z' = e^d x + i e^{-d} y$ defining an isometric embedding $[0,\infinity) \to \Teich$. The image of this map is called the \emph{\Teichmuller\ ray with basepoint $\sigma$ and direction $P[\fol]$}.

A measured foliation $\fol$ is \emph{arational}\index{arational measured foliation} if $\<\fol,c\> > 0$ for every $c \in \C$. Equivalently, there is no closed loop in $S$ which is everywhere parallel to $\fol$, nor is there any closed path from puncture to puncture which is everywhere parallel to $\fol$. If $\fol,\fol'$ are arational, and if $\<\fol,\fol'\> > 0$, then $\fol$, $\fol'$ jointly fill $S$.

\paragraph{Virtual torsion free.} Serre \cite{Serre:rigidite} originally proved that the group $\MCG$ has a finite index subgroup which is torsion free. For another proof see \cite{Ivanov:subgroups} Corollary~1.5.

\begin{theorem}
\label{TheoremVTF}
The kernel of the homomorphism $\MCG(S) \to \GL(H_1(S;\Z/3))$ is torsion free. This kernel has index bounded above by $\abs{\GL(H_1(S;\Z/3)))}$ which is finite because $H_1(S;\Z/3)$ has finite rank $\le 2g+p$ over $\Z/3$.
\end{theorem}

Note that the bound in the index of the kernel is exponential in $g$ and $p$.

\paragraph{Gromov hyperbolic metric spaces.} 
Much of this survey is motivated by analogues in the theory of discrete groups of isometries of Gromov hyperbolic metric spaces, particularly hyperbolic spaces $\hyp^n$. We recall some of the notions of this theory, whose details can be found in, for example, \cite{Gromov:hyperbolic}, \cite{GhysHarpe:hyperboliques}, \cite{CDP}, \cite{Cannon:TheoryHyp}.

Let $X$ be a metric space. $X$ is \emph{proper} if closed balls are compact. $X$ is a \emph{geodesic metric space} if for any $x,y \in X$ there exists an isometrically embedded path $\gamma \from [0,d(x,y)] \to X$, a \emph{geodesic}, from $x$ to $y$. We'll use the notation $[x,y]$ to refer to the image of some geodesic. The \emph{Hausdorff distance} between two sets $A,B \subset X$ is the infimum of $\delta>0$ such that $A \subset N_\delta(B)$ and $B \subset N_\delta(A)$.

A geodesic metric space $X$ is \emph{hyperbolic}\index{hyperbolic metric space} in the sense of Gromov if there exists $\delta>0$ such that for any $x,y,z$ and any geodesics $[x,y]$, $[y,z]$, $[z,x]$ in $X$, $[x,y]$ is contained in the $\delta$-neighborhood of $[y,z] \union [z,x]$. The \emph{boundary} $\bdy X$ is defined to be the set of geodesic rays modulo the relation of finite Hausdorff distance; given a geodesic ray $\rho$, let $[\rho]$ denote the corresponding point of $\bdy X$. There is a natural \emph{Gromov topology} on $X \union \bdy X$, whose definition we shall not recount, but we shall mention several important properties:
\begin{itemize}
\item If $X$ is a proper geodesic metric space then $\bdy X$ and $X \union \bdy X$ are compact.
\item Every geodesic ray $\rho$ converges in $X \union \bdy X$ to $[\rho]$.
\item If $X$ is proper then for any $x \in X$ and any $\xi \in \bdy X$ there exists a geodesic ray $\rho$ based at $x$ and converging to $\xi$, and for any $\xi,\eta \in \bdy X$ there exists a bi-infinite geodesic whose two ends converge to $\xi,\eta$.
\end{itemize}

Of course, the primary example of the Gromov topology is the union of $\hyp^n$ with its sphere at infinity $S^{n-1}_\infinity$. When $X$ is not proper then the compactness properties above can fail. 

Another important example for us of a Gromov hyperbolic metric space is the curve complex of $S$, discussed below.

\section{The classification of mapping classes}
\label{SectionClassification}

The discrete cyclic groups of isometries of $\hyp^n$ or of any proper, geodesic, Gromov hyperbolic metric space $X$ fall into three classes: \emph{elliptic} groups which are finite and which fix a point of $X$; \emph{parabolic} groups which are infinite cyclic and which fix a point of $\bdy X$ but no point of $X$; and \emph{loxodromic} groups which are infinite cyclic and which translate along a bi-infinite geodesic of $X$. 

Thurston discovered an analogous classification of elements of $\MCG$ \cite{Thurston:surfaces}; see also \cite{FLP} and \cite{CassonBleiler}. Work of Bers strengthens the analogy in terms of the action on $\Teich$ \cite{Bers:ThurstonTheorem}.

There are two major classifications of a mapping class $\phi \in \MCG$: whether $\phi$ is finite order; and whether $\phi$ is \emph{reducible}\index{reducible mapping class} meaning that there exists an essential curve family $\C$, called a \emph{reduction family} for $\phi$, such that $\phi[\C]=[\C]$. If $\phi$ is not reducible then it is \emph{irreducible}. 

A homeomorphism $\Phi \from S \to S$ is said to be \emph{pseudo-Anosov}\index{pseudo-Anosov homeomorphism} if there exists a transverse pair of measured foliations $\fol^s$, $\fol^u$ which jointly fill the surface, called the \emph{stable and unstable measured foliations} of $\Phi$, and there exists a number $\lambda>1$ called the \emph{expansion factor} of $\Phi$, such that $\Phi$ takes each leaf of $\fol^s$ to a leaf of $\fol^s$, each leaf of $\fol^u$ to a leaf of $\fol^u$, $\Phi$ contracts the leaves of $\fol^s$ by a factor of $\lambda$ with respect to the $\fol^u$ measure, and $\Phi$ expands the leaves of $\fol^u$ by a factor of $\lambda$ with respect to the $\fol^s$ measure. To say this more briefly, $\Phi(\fol^s)=\lambda\fol^s$ and $\Phi(\fol^u)=\lambda^\inv \fol^u$. 

A mapping class $\phi \in \MCG(S)$ is said to be pseudo-Anosov if it is represented by a pseudo-Anosov homeomorphism.

Thurston proved that every mapping class $\phi \in \MCG(S)$ falls into one of three types: finite order, reducible, or pseudo-Anosov; in other words, every infinite order, irreducible mapping class is \emph{pseudo-Anosov}. The first two types are not mutually exclusive, but a mapping class is pseudo-Anosov if and only if it is infinite order and irreducible.

There are several necessary and/or sufficient conditions for the three types of mapping classes, expressed in terms of the geometry and dynamics of the action on compactified \Teichmuller\ space.

\paragraph{Finite order mapping classes.} A mapping class $\phi \in \MCG(S)$ is of finite order if and only if it has a fixed point in $\Teich$, in which case there exists a conformal structure $\sigma$ on $S$ and a homeomorphism $\Phi$ representing $\phi$ such that $\Phi(\sigma)=\sigma$. It follows that $\Phi$ is a finite order homeomorphism.

A sufficient condition for $\phi \in \MCG(S)$ to be of finite order, proved in \cite{FLP}, is the existence of an arational measured foliation $\fol$ such that $\phi[\fol]=[\fol]$.
 
It is possible for a finite order mapping class $\phi$ to be irreducible. This happens if and only if, for $\Phi$ a finite order homeomorphism representing $\phi$, the quotient orbifold $S/\Phi$ is a sphere with three cone points.

There is a bound on the order of a finite order element of $\MCG(S)$, in fact a bound on the order of a finite subgroup $G \subgroup \MCG(S)$, depending only on the topology of $S$. One way to see this is by applying Theorem~\ref{TheoremVTF}, which gives a bound which is exponential in $g$ and $p$. A linear bound is obtained using the much deeper theorem of Gabai \cite{Gabai:convergence} and Casson--Jungreis \cite{CassonJungreis} which says that $G$ is realized as a finite group of homeomorphisms of $S$, also denoted $G$. We therefore have $\abs{G} = \chi(S) / \chi(S/G) = (2-2g-p)/\chi(S/G)$, where the denominator is the orbifold Euler characteristic of the quotient orbifold $S/G$. The numerator and denominator of this fraction being negative, $\abs{G}$ is maximized when $\chi(S/G)$ is maximized.

The maximum value of the Euler characteristic of a compact, oriented hyperbolic \nb{2}orbifold is realized by the spherical orbifold with three cone points of angles $\frac{2\pi}{2}, \, \frac{2\pi}{3}, \, \frac{2\pi}{7}$, whose Euler characteristic equals $2 - (1-\frac{1}{2}) - (1 - \frac{1}{3}) - (1 - \frac{1}{7}) = \frac{1}{42}$. It follows that $\abs{G} \le 42(2g-2)$ when $S$ is compact.

When $S$ has one or more punctures there is a slightly better bound, because the maximum value of the Euler characteristic of a noncompact, finite area, oriented hyperbolic \nb{2}orbifold is acheived by the modular space of $\SL(2,\Z)$, a once-punctured sphere with two cone points of angles $\frac{2\pi}{2}, \, \frac{2\pi}{3}$ and with Euler characteristic $2 - (1 - \frac{1}{2}) - (1 - \frac{1}{3}) = \frac{1}{6}$, and so $\abs{G} \le 6(2g-2+p)$.

To summarize:
\begin{theorem}
\label{TheoremFiniteBound}
If $G$ is a finite subgroup of $\MCG$ then
$$\abs{G} \le \begin{cases}
84g-84 &\text{if}\quad p=0 \\
12g-12+6p &\text{if}\quad p>0
\end{cases}
$$
\qed\end{theorem}

\paragraph{Pseudo-Anosov mapping classes.} Bers proved \cite{Bers:ThurstonTheorem} that a mapping class $\phi \in \MCG(S)$ is pseudo-Anosov if and only if the infimum of $d([\sigma],\phi[\sigma])$ is positive and is acheived by some $[\sigma] \in \Teich$. 
In this case, by combining with the fact that the stable and unstable foliations $\fol^s,\fol^u$ are uniquely ergodic, it follows that the set of points $[\sigma]$ where the infimum is acheived is a bi-infinite geodesic in $\Teich$ along which $\phi$ translates, with translation distance equal to $\log(\lambda)$. This geodesic is called the \emph{axis} for $\phi$ in $\Teich$. This axis has two ideal endpoints in $\PMF$, namely $P[\fol^s]$ and $P[\fol^u]$. The action of $\phi$ on $\overline\Teich$ is a \emph{source sink action}, with \emph{source} $P[\fol^s]$ and \emph{sink} $P[\fol^s]$, meaning that for every point $x \in \overline\Teich - \{P[\fol^s],P[\fol^u]\}$, 
$$\lim_{i \to -\infinity} \phi^i(x) = P[\fol^s] \quad\text{and}\quad \lim_{i \to +\infinity} \phi^i(x) = P[\fol^u]
$$
It follows that the fixed point set of $\phi$ in $\overline\Teich$ is
$$\Fix(\phi) = \{P[\fol^s],P[\fol^u]\}
$$

A sufficient (and necessary) condition for $\phi \in \MCG(S)$ to be pseudo-Anosov is the existence of an arational measured foliation $\fol$ and $\lambda>0$, $\lambda\ne 1$ such that $\phi[\fol]=\lambda[\fol]$; if this happens then $P[\fol] = P[\fol^s]$ or $P[\fol^u]$. See \cite{FLP} for the proof.

\paragraph{Reducible mapping classes.} If a mapping class $\phi$ is reducible, and if $\C$ is a reducing family for $\phi$, then there exists a representative $\Phi$ of $\phi$ such that $\Phi(\C)=\C$. Given a component $S'$ of $S-\C$, let $n \ge 1$ be the first return time of $\Phi$ to $S'$, meaning that $\Phi^n(S')=S'$ and $\Phi^i(S') \ne S'$ for $0<i<n$. Then the first return time $n$ and the mapping class of the first return map $\Phi^n$ are determined by $\phi$, independent of the representative $\Phi$. These mapping classes are called the \emph{component mapping classes} of $\phi$ relative to $\C$.

A reducing family $\C$ for a mapping class $\phi$ is \emph{complete} if all of the component mapping classes of $\phi$ relative to $\C$ are finite order or pseudo-Anosov. A minimal, complete reducing family $\C$ for $\phi$ is unique up to isotopy; see for example \cite{BirmanLubotzkyMcCarthy}, \cite{HandelThurston:nielsen}. If $\C$, $\C'$ are any two minimal, complete reducing families for $\phi$, if $\Phi,\Phi'$ are representatives of $\phi$ such that $\Phi(\C)=\C$, $\Phi'(\C')=\C'$, and if $\Psi \in \Homeo_0(S)$ satisfies $\Psi(\C)=\C'$, then $\Psi$ conjugates the component mapping classes of $\phi$ relative to $\C$ to the component mapping classes of $\phi$ relative to $\C'$. We may therefore speak, in a well-defined manner, of the \emph{component mapping classes of $\phi$}, meaning the component mapping classes of $\phi$ relative to any minimal, complete reducing family. 

Bers proved \cite{Bers:ThurstonTheorem} that a mapping class $\phi$ is infinite order and reducible if and only if the infimum of $d([\sigma],\phi[\sigma])$ is not acheived by any $[\sigma] \in \Teich$. Moreover, two subcases are distinguished: if the infimum equals zero then each component mapping class of $\phi$ has finite order; whereas if the infimum is positive then some component mapping class of $\phi$ is pseudo-Anosov.

A sufficient (and necessary) condition for $\phi \in \MCG(S)$ to be reducible is the existence of a non-arational measured foliation $\fol$ such that $\phi[\fol]=\fol$.

A complete description of the dynamics of a reducible mapping class acting on $\PMF = \bdy\Teich$ is obtained in \cite{Ivanov:subgroups}. 

\paragraph{Independence of pseudo-Anosov homeomorphisms.} In a discrete group acting on $\hyp^n$, given two loxodromic elements $\phi,\phi'$, either $\Fix(\phi)=\Fix(\phi')$ or $\Fix(\phi) \intersect \Fix(\phi') = \emptyset$. The same is true for pseudo-Anosov homeomorphisms. Proofs are given in \cite{Ivanov:subgroups,McCarthyPapa:dynamics}. 
We shall give a proof here.

\begin{lemma}
\label{LemmaDisjointFixed}
For any pseudo-Anosov $\phi,\phi' \in \MCG(S)$, either $\Fix(\phi)=\Fix(\phi')$ or $\Fix(\phi) \intersect \Fix(\phi') = \emptyset$. In the latter case we say that $\phi,\phi'$ are \emph{independent}.
\end{lemma}

\begin{proof} Supposing that $\Fix(\phi) \intersect \Fix(\phi') \ne \emptyset$, up to inverses of $\phi,\phi'$ and up to isotopy we may assume that $\fol := \fol^u(\phi) = \fol^u(\phi') $, and that $\phi[\fol] = \lambda[\fol]$ and $\phi'[\fol]=\lambda'[\fol]$ for some $\lambda,\lambda' > 1$. Let $\Gamma$ be the intersection of $\Stab(P[\fol])$ and the kernel of the homomorphism $\MCG(S) \to \GL(H_1(S;\Z/3))$, and so by Theorem~\ref{TheoremVTF} it follows that $\Gamma$ has finite index in $\Stab(P[\fol])$. It suffices to prove that $\Gamma$ is infinite cyclic, for in that case $\phi,\phi'$ have positive powers in $\Gamma$, implying that $\phi^m=\phi'{}^n$ for some $m,n \ge 0$, implying that $[\fol^s(\phi)] = [\fol^s(\phi^m)] = [\fol^s(\phi'{}^n)] = [\fol^s(\phi')]$.

Consider the homomorphism $\ell \from \Gamma \to \reals_+$ which is characterized by the equation $\psi[\fol] = \ell(\psi)[\fol]$. Suppose that $\psi \in \kernel(\ell)$. As seen above, $\psi$ has finite order. But $\psi \in \kernel(\MCG(S) \to \GL(H_1(S;\Z/3)))$ and so by Theorem~\ref{TheoremVTF}, $\psi$ is the identity. The homomorphism $\ell$ is thus injective. 
But the image of $\ell$ is discrete --- indeed, the set of expansion factors of pseudo-Anosov elements of $\MCG(S)$ is bounded away from $1$, because each is an eigenvalue of a square integer matrix of bounded size \cite{FLP}. It follows that $\Gamma$ is infinite cyclic.
\end{proof}

\paragraph{The ping-pong argument and free subgroups.} Here is the central construction of the Tits alternative:

\begin{lemma}
\label{LemmaPsAnPingPong}
For any two independent pseudo-Anosov mapping classes $\phi,\psi$, for sufficiently large integers $m,n>0$ the elements $\phi^m,\psi^n$ freely generate a free subgroup of $\MCG(S)$.
\end{lemma}

The proof of this lemma is a standard ping-pong argument, using that $\phi,\psi$ each act with source-sink dynamics on $\Teich$, and that their fixed points sets are disjoint. A general ping-pong lemma stated in a topological context is as follows:

\begin{lemma}[Topological Ping Pong]
\index{ping-pong!topological}
Let $G$ act on a topological space $X$, let $\phi_1,\ldots,\phi_k \in G$, and suppose that there exists a pairwise disjoint collection of open sets $U^-_i, U^+_i$, $i=1,\ldots,k$ such that for each $i \in \{1,\ldots,k\}$ and each $\epsilon \in \{-,+\}$ we have
$$\phi^\epsilon_i(X-U^{-\epsilon}_i) \subset U^\epsilon_i
$$
Then $G$ is freely generated by $\phi_1,\ldots,\phi_k$.
\end{lemma}

The way that Lemma~\ref{LemmaPsAnPingPong} is proved is to note that when $\phi,\psi$ are independent pseudo-Anosov mapping classes, and if $U^s_\phi$, $U^u_\phi$, $U^s_\psi$, $U^u_\psi$ are pairwise disjoint neighborhoods of the fixed points $P[\fol^s_\phi]$, $P[\fol^u_\phi]$, $P[\fol^s_\psi]$, $P[\fol^u_\psi]$, then by applying the source--sink dynamics of $\phi$ and of $\psi$, the hypotheses of the Ping Pong Lemma are satisfied for $\phi^m,\psi^n$ as long as $m,n$ are sufficiently large.

One of the main themes of the rest of this survey is that by taking still more care with the proof of Lemma~\ref{LemmaPsAnPingPong}, one can show that the free subgroup $\<\phi^m,\psi^n\>$ satisfies stronger geometric properties reminiscent of convex cocompact discrete subgroups of hyperbolic isometries.

\paragraph{The virtual centralizer and virtual normalizer of a pseudo-Anosov cyclic subgroup.} Consider a group $G$ and a subgroup $H \subgroup G$. The \emph{virtual centralizer} $VC(H)$ of $H$ in $G$ is the subgroup of all $g \in G$ which commute with a finite index subgroup of $H$. The \emph{virtual normalizer} $VN(H)$ is the subgroup of all $g \in G$ such that $gHg^\inv \intersect H$ has finite index in $gHg^\inv$ and in $H$.

The virtual centralizer and virtual normalizer of an infinite cyclic pseudo-Anosov subgroup each have a particularly nice geometric/dynamic description, and the proof of this description uses the ping-pong argument in a nice way. Given an action of a group $G$ on a space $X$, for a subspace $S \in X$ denote $\Stab(S) = \{g \in G \suchthat g(S)=S\}$. If $S$ is homeomorphic to a line denote $\Stab_+(S)$ to be the subgroup of $\Stab(S)$ that preserves the orientation on $S$.

\begin{theorem}
\label{TheoremVCVN}
Given a pseudo-Anosov $\phi \in \MCG(S)$ with stable and unstable foliations $[\fol^s]$, $[\fol^u]$ and \Teichmuller\ axis $A$, we have:
\begin{align*}
VC\<\phi\> &= \Stab(P[\fol^s]) = \Stab(P[\fol^u]) = \Stab_+(A) \\
VN\<\phi\> &= \Stab\{P[\fol^s],P[\fol^u]\} = \Stab(A)
\end{align*}
The group $\Stab_+(A)$ 
has index at most $2$ in $\Stab(A)$, it has an infinite cyclic subgroup of finite index, and the index in $\Stab_+(A)$ of the maximal such subgroup is bounded by the maximum order of a finite subgroup of $\MCG(S)$.
\end{theorem}

\begin{proof}
The group $\Stab(A)$ acts properly discontinuously on the \Teichmuller\ geodesic $A$ and so is virtually cyclic. The inclusions $\Stab_+(A) \subset VC\<\phi\>$ and $\Stab(A) \subset VN\<\phi\>$ follow immediately. The inclusions $\Stab(P[\fol^s]), \Stab(P[\fol^u]) \subset \Stab_+(A)$ and the equation $\Stab\{P[\fol^s],P[\fol^u]\} = \Stab(A)$ are obvious. 

Suppose that $\psi \in \MCG(S) - \Stab\{P[\fol^s],P[\fol^u]\}$. Then $\phi'=\psi\phi\psi^\inv$ is pseudo-Anosov, and it is independent of $\phi$ by Lemma~\ref{LemmaDisjointFixed}. Some powers of $\phi$ and $\phi'$ therefore generate a free group of rank~2, and hence $\psi \not\in VN\<\phi\>$.

Suppose next that $\psi \in \MCG(S) -  \Stab(P[\fol^u])$, and again consider the pseudo-Anosov mapping class $\phi'=\psi\phi\psi^\inv$. If $\psi(P[\fol^u]) \ne P[\fol^s]$ then $\phi'$ independent of $\phi$ and $\psi \not\in VN\<\phi\>$, so $\psi\not\in VC\<\phi\>$. 
If $\psi(P[\fol^u]) = P[\fol^s]$ but $\psi(P[\fol^s]) \ne P[\fol^u]$ then $\Fix(\phi)$, $\Fix(\phi')$ are not disjoint and not equal, violating Lemma~\ref{LemmaDisjointFixed}. And if $\psi$ interchanges $P[\fol^u]$ and $P[\fol^s]$ then $\phi'$ translates along $A$ in the direction opposite to $\phi$ and so $\psi \not\in VC\<\phi\>$.
\end{proof}

\section{The Tits alternative and abelian subgroups}
\label{SectionTits}

The following theorem was proved independently by Ivanov \cite{Ivanov:subgroups} and McCarthy \cite{McCarthy:Tits}:

\begin{theorem}[The Tits Alternative]
\index{Tits alternative}
For each finite type surface $S$ there exist constants $I,R$ such that for any subgroup $G$ of $\MCG(S)$, either $G$ contains a free subgroup of rank~$2$, or $G$ contains an abelian subgroup of index $\le I$ and rank $\le R$.
\end{theorem}

Ivanov gave an underlying geometric description of each subgroup of $\MCG(S)$, from which the above Tits alternative follows quickly. A subgroup $G \subgroup \MCG(S)$ is \emph{reducible} if there exists an essential curve family which is invariant under the action of each element of $G$. If no such curve family exists then $G$ is \emph{irreducible}.

\begin{theorem}
\label{TheoremIvanov}
If $G$ is an infinite, irreducible subgroup of $\MCG(S)$ then $G$ contains a pseudo-Anosov element.
\end{theorem}

For the proof see Ivanov's book \cite{Ivanov:subgroups}.

\begin{proof}[Proof of Tits alternative, from Theorem \ref{TheoremIvanov}]
Let $G$ be an arbitrary subgroup of $\MCG$.

Suppose first that $G$ is irreducible. Applying  Theorem~\ref{TheoremIvanov}, $G$ contains a pseudo-Anosov mapping class $\phi$. If $G \subset \Stab(\Fix(\phi)) = VN(\phi)$ then by Theorem~\ref{TheoremVCVN} the group $G$ has an infinite cyclic subgroup of index at most twice the maximum size of a finite subgroup of $\MCG$. If $G\not\subset\Stab(\Fix(\phi))$, say $\psi \in G - \Stab(\Fix(\phi))$, then by Lemma~\ref{LemmaDisjointFixed} $\psi\phi\psi^\inv$ is a pseudo-Anosov element independent of $\phi$ and so by the ping pong argument $G$ contains a free group of rank~2.

Suppose now that $G$ is reducible. Let $\C$ be a maximal essential curve family invariant under $G$. The group $G$ acts on the set of oriented elements of $\C$; let $G_\C$ denote the kernel of this action, a finite index subgroup of $G$ whose index is bounded above by $2^n\abs{\C}!$. Letting $\{S_i\}$ be the components of $S-\C$, it follows that $G_\C$ preserves each $S_i$ up to isotopy, and so there is a well-defined restriction homomorphism $r_i \from G_\C \to \MCG(S_i)$ whose image we denote $G \restrict S_i$. By maximality of $\C$ it follows that $G \restrict S_i$ is irreducible. If $G \restrict S_i$ has two independent pseudo-Anosov elements then $G \restrict S_i$ has a free subgroup of rank $\ge 2$, and so therefore does $G$. 

We have reduced to the case that $G \restrict S_i$ is either finite or is contained in $VN(\phi_i)$ for some pseudo-Anosov element $\phi_i \in \MCG(S_i)$. In this case we find a free abelian subgroup of finite index in $G$. Let $A_i \subset G \restrict S_i$ denote either the trivial subgroup or an infinite cyclic group, with index bounded as in Theorem~\ref{TheoremFiniteBound} or~\ref{TheoremVCVN}. Consider the subgroup $H = \cap_i r^\inv(A_i) \subgroup G$. The number of subsurfaces $S_i$ is bounded in terms of $g$ and $p$, and so we obtain a finite bound on the index of $H$ in $G$, in terms of $g$ and $p$. Restricted to each $S_i$, the subgroup $H$ is either trivial or infinite cyclic, in the latter case generated by a pseudo-Anosov element of $S_i$. $H$ also contains Dehn twists about the curves in $\C$. These pseudo-Anosov subsurface mapping classes and the Dehn twists all commute with each other. Thus $H$ is a free abelian group, whose rank is bounded by $\abs{\C} + \abs{\{S_i\}}$, which is bounded in terms of $g$ and $p$.
\end{proof}

Note that this proof gives a complete description of all free abelian subgroups of $\MCG(S)$: each irreducible one is infinite cyclic generated by a pseudo-Anosov mapping class; and each reducible one is freely generated by Dehn twists along curves $\C$ and pseudo-Anosov mapping classes on components $\{S_i\}$ of $S-\C$. To compute the maximal rank exactly, note that on any subsurface $S_i$ that supports a pseudo-Anosov homeomorphism, one can replace that homeomorphism by a Dehn twist, enlarging $\C$ without decreasing the rank. The maximum rank is thus acheived when each $S_i$ is a pair of pants, in which case the rank equals the maximum value of $\abs{\C}$, which equals $3g-3+p$.

\section{Limit sets}
\label{SectionLimits}

\index{limit set}
Given a Gromov hyperbolic, geodesic metric space $X$ and a finitely generated, discrete subgroup $G \subgroup \Isom(X)$, its \emph{limit set} $\Lambda_G \subset \bdy X$ is, by definition, the accumulation set in $\bdy X$ of any orbit of $G$ in $X$. As long as $G$ is nonelementary (that is, not virtually abelian), the limit set may be characterized as the closure of the set of fixed points of loxodromic elements of $G$, or as the unique minimal closed nonempty $G$-invariant subset of $\bdy X$. Assuming in addition that $X$ is proper, it follows that $\Lambda_G$ is compact, and the \emph{domain of discontinuity} $\Delta_G = \bdy X - \Lambda_G$ is characterized as the unique maximal open $G$-invariant set on which the action of $G$ is properly discontinuous.

In \cite{Masur:handlebodies}, Masur investigated notions of the limit set and domain of discontinuity for the \emph{handlebody subgroup} of $\MCG(S)$ when $S$ is the boundary of a compact genus $g$ handlebody $H$; this subgroup consists of the mapping classes of homeomorphisms of $S$ that extend to $H$.

In \cite{McCarthyPapa:dynamics}, McCarthy and Papadopoulos generalized Masur's results to arbitrary subgroups of $\MCG(S)$. The results are not quite as pretty as they are for discrete groups of hyperbolic isometries. But as we shall see, the results of Kent and Leininger \cite{KentLeininger:Shadows} show that limit sets of convex cocompact subgroups of $\MCG$ behavior very much like their counterparts in hyperbolic geometry.

We start with some elementary cases. When $G$ is a finite subgroup of $\MCG$, the natural choice of the limit set is $\Lambda_G=\emptyset$, and the action of $G$ on $\Delta_G=\PMF$ is properly discontinuous. When $G$ is virtually cyclic, containing a pseudo-Anosov cyclic subgroup $\<\phi\>$ with finite index, let $\Lambda_G = \{P[\fol^s_\phi],P[\fol^u_\phi]\}$. In the case where the action of $G$ does not fix each element of $\Lambda_G$, that is when $VC\<\phi\>$ is an index~2 subgroup of $VN\<\phi\>$, then $\Lambda_G$ is the unique nonempty minimal closed $G$-invariant subset of $\PMF$. In the other case, when $VC\<\phi\>=VN\<\phi\>$, then the two points of $\Lambda_G$ are the only nonempty minimal closed $G$-invariant subsets of $\PMF$. In either case, the action of $G$ on $\Delta_G = \PMF-\Lambda_G$ is properly discontinuous.

Consider now the case where $G$ is an infinite, irreducible subgroup of $\MCG$ and $G$ has no finite index cyclic subgroup. Define $\Lambda_0 \subset \PMF$ to be the set of fixed points of pseudo-Anosov elements of $G$. Define $\Lambda_G \subset \PMF$ to be the closure of $\Lambda_0$. Define $Z\Lambda_G$ to be the ``zero set'' of $\Lambda_G$, the set of all $P[\fol] \in \PMF$ such that for some $P[\fol'] \in \Lambda_G$ we have 
$\<[\fol],[\fol']\> = 0$; this is well-defined independent of the choice of $[\fol],[\fol']$ in their projective Whitehead classes. Define $\Delta_G = \PMF - Z\Lambda_G$.

\begin{theorem}[\cite{McCarthyPapa:dynamics}]
\label{TheoremMcPapaIrred}
For $G$ an infinite, irreducible subgroup of $\MCG$, the set $\Lambda_G$ is the unique minimal nonempty closed $G$-invariant subset of $\PMF$. If $\Lambda_G \ne \PMF$ then $\Delta_G \ne \emptyset$, and the action of $G$ on $\Delta_G$ is properly discontinuous.
\end{theorem}

We turn now to the case that $G$ is an infinite, reducible subgroup of $\MCG(S)$. In general there will be no unique closed minimal $G$-invariant subset of $\PMF$. For example, if $G$ is a free abelian group generated by Dehn twists about the curves in an essential curve family $\C$, then $G$ has many, many fixed points, namely any point in $\PMF$ having zero intersection number with each curve in $\C$.  Nevertheless, as McCarthy and Papadopoulos show, there is still a reasonable candidate for a limit set and a domain of discontinuity.

Let $\C$ be an essential curve family invariant under $G$ such that if $\{S_i\}$ is the set of components of $S-\C$, then the restriction $G \restrict S_i \subgroup \MCG(S_i)$ is an irreducible subgroup of $\MCG(S_i)$. As is shown in, say, \cite{McCarthyPapa:dynamics}, there exists a unique such family $\C$ which is minimal with respect to inclusion. Reindex the $S_i$ so that $S_1,\ldots,S_n$ are the surfaces on which the restriction groups $G \restrict S_i$ are infinite. Each of the groups $G \restrict S_i$ for $i=1,\ldots,n$ therefore contains a pseudo-Anosov element of $\MCG(S_i)$. The stable and unstable foliations of pseudo-Anosov elements of $G \restrict S_i$, being measured foliations on $S_i$, may by enlargement be regarded as measured foliations on $S$; let $\Lambda^i_0 \subset \PMF$ be the set of all points in $\PMF$ obtained in this manner. Let $\Lambda^i \subset \PMF$ be the closure of $\Lambda^i_0$. Define
$$\Lambda_G = \left(\bigcup_{i=1}^n \Lambda^i\right) \union \left(\bigcup_{c \in \C} P[c] \right)
$$
Define the zero set $Z\Lambda_G$ and the domain $\Delta_G = \PMF-Z\Lambda_G$ as above.

\begin{theorem}[\cite{McCarthyPapa:dynamics}]
For each infinite, reducible subgroup of $\MCG$, the set $\Lambda_G$ is a nonempty closed $G$-invariant subset of $\PMF$. If $\Lambda_G \ne \PMF$ then $\Delta_G \ne \emptyset$, and the action of $G$ on $\Delta_G$ is properly discontinuous.
\end{theorem}

\index{domain of discontinuity}
The set $\Lambda_G$ shall be called the \emph{limit set} and $\Delta_G$ the \emph{domain of discontinuity} of $G$. 

\section{Convex cocompact subgroups}
\label{SectionConvexCocompact}
\index{convex cocompact subgroup}
Given a finitely generated, discrete subgroup $G \subgroup \Isom(\hyp^n)$, the \emph{convex hull} of its limit set $\Lambda_G$ is the smallest closed, convex subset $\H_G \subset \hyp^n$ whose set of accumulation points in $S^{n-1}_\infinity$ is $\Lambda_G$. The action of $G$ on $\hyp^n$ preserves $\H_G$, and $G$ is said to be \emph{convex cocompact} if $G$ acts cocompactly on $\H_G$. Several conditions are equivalent to convex cocompactness: some (every) orbit of $G$ in $\H_G$ is quasiconvex; the action of $G$ on $\Delta_G$ is cocompact; the action of $G$ on $\hyp^n \union \Delta_G$ is cocompact. Moreover, if $G$ is convex cocompact then $G$ is word hyperbolic, there is a $G$-equivariant homeomorphism from the Gromov boundary $\bdy G$ to $\Lambda_G$, and this homeomorphism extends to a $G$-equivariant continuous map from the Gromov compactification $G \union \bdy G$ to $\H_G \union \Lambda_G$.

These notions can be generalized with care to a finitely generated, discrete group $G$ of isometries of any Gromov hyperbolic, geodesic metric space~$X$. 

The theory of convex cocompact subgroups of $\MCG$, originated in \cite{FarbMosher:quasiconvex} and further developed in \cite{KentLeininger:Shadows}, and \cite{Hamenstadt:extensions}, gives analogues to these results about convex cocompact subgroups of $\MCG$. In addition, the theory develops further equivalent characterizations of convex compactness that are special to the setting of $\MCG$ and are important for applications. Thus far, the only examples of convex cocompact subgroups of $\MCG$ are free \cite{FarbMosher:quasiconvex} or virtually free examples as constructed by Honglin Min in her dissertation \cite{Min:HyperbolicGraphs}.

We define a subset $X \subset \Teich$ to be \emph{quasiconvex} if there exist a constant $C$ such that for every $x,y \in X$, each point of the \Teichmuller\ geodesic $[x,y]$ is within distance $C$ of some point of $X$.

Given a closed subset $\Lambda \subset \PMF$, if every pair $\xi \ne \eta \in \Lambda$ jointly fills $S$ then we define the \emph{weak hull} of $\Lambda$ to be the subset of $\Teich$ defined by
$$\H(\Lambda) = \union\{\geodesic{\xi}{\eta} \suchthat \xi \ne \eta \in \Lambda\}
$$
and we say that the \emph{weak hull of $\Lambda$ is defined}. Otherwise, if there exist $\xi \ne \eta \in \Lambda$ which do not jointly fill $S$ then we say that the weak hull of $\Lambda$ is undefined.

\begin{theorem}[\cite{FarbMosher:quasiconvex}]
For any finitely generated subgroup $G \subgroup \MCG$, the following are equivalent.
\begin{itemize}
\item Every (some) orbit of $G$ on $\MCG$ is quasiconvex in $\Teich$.
\item $G$ is word hyperbolic, and there exists a $G$-equivariant embedding $\bdy G \hookrightarrow \PMF$ with image denoted $\Lambda$, such that the weak hull of $\Lambda$ is defined, the action of $G$ on $\Hull(\Lambda)$ is cocompact, and the extension of the map $\bdy G \hookrightarrow \PMF$ by any $G$-equivariant map $G \to \Hull(\Lambda)$ is a continuous map $G \union \bdy G \to \Hull(\Lambda) \union \Lambda$.
\end{itemize}
\end{theorem}

A group satisfying these conditions is called \emph{convex cocompact}. When $G$ is convex cocompact then every infinite order element of $G$ is pseudo-Anosov.

The connection of convex cocompactness with the limit set $\Lambda_G$ and domain of discontinuity $\Delta_G$ as defined by McCarthy-Papadapoulos is given in the following results of Kent and Leininger: 

\begin{theorem}[\cite{KentLeininger:Shadows}]
For any finitely generated subgroup $G \subgroup \MCG$, the following are equivalent:
\begin{itemize}
\item $G$ is convex cocompact.
\item The weak hull $\Hull_G$ of the McCarthy-Papadapoulos limit set $\Lambda_G$ is defined, and the action of $G$ on $\Hull_G$ is cocompact.
\item The action of $G$ on the McCarthy-Papadapoulos domain of discontinuity $\Delta_G$ is cocompact.
\item The action of $G$ on $\Teich \union \Delta_G$ is cocompact.
\end{itemize}
Furthermore, if $G$ is convex cocompact, then the image of the embedding $\bdy G \to \PMF$ is the McCarthy-Papadapoulos limit set $\Lambda_G$, and $Z\Lambda_G = \Lambda_G$, and so $\Delta_G = \PMF - \Lambda_G$.
\end{theorem}

\paragraph{Convex cocompactness on the curve complex.} \index{curve complex}
The \emph{curve complex} $\CC=\CC(S)$ is the simplicial complex whose vertex set is $\C$,  and whose $n$-simplices are the isotopy classes of essential curve families of cardinality $n+1$. To be precise, a set of distinct vertices $[c_0],\ldots,[c_n] \in \C$ spans an $n$-simplex of $\CC$ if and only if there are representatives $c_0,\ldots,c_n$ which are pairwise disjoint. The complex $\CC$ is locally infinite, in fact the link of every simplex of positive codimension is infinite. 

The natural action of $\MCG$ on $\CC$ has finitely many orbits of simplices. The number of orbits of vertices is
$$
m = \begin{cases}
1 + \ceiling{\frac{g-1}{2}} &\text{if}\quad g \ge 1, p=0 \\
\ceiling{\frac{p-3}{2}}       &\text{if}\quad g=0, p \ge 4\\
1 + \ceiling{\frac{(g+1)(p+1)-4}{2}}  &\text{if}\quad g,p \ge 1
\end{cases}
$$
When $g \ge 1$ there is one orbit of nonseparating curves, but in genus zero every curve is separating. The orbit of a separating curve is determined by how it partitions the genus and the number of punctures. When $p=0$, the relevant partition is $g=g_1+g_2$ with $g_1,g_2 \in \{1,\ldots,g\}$, and when $g=0$ the relevant partition $p=p_1+p_2$ with $p_1,p_2 \in \{2,\ldots,p\}$. When $g \ge 1$, $p \ge 1$ then the relevant partition is an unordered partition of the ordered pair $(g,p) = (g_1,p_1) + (g_2,p_2)$, where the four values $(g_i,p_i) = (0,0)$, $(g,p)$, $(0,1)$, $(g,p-1)$ are forbidden but all other values with $g_i \in \{0,\ldots,g\}$ and $p_i \in \{0,\ldots,p\}$ are allowed. 

The large scale geometric properties of $\CC$ are of $\MCG$ are related in the following manner. 

Given a finitely generated group $G$ and a finite collection of finitely generated subgroups $H_1,\ldots,H_m$, assume that the generators of $G$ include generators of each of the subgroups $H_1,\ldots,H_m$, so that the Cayley graph of $H_i$ embeds, $H_i$-equivariantly, in the Cayley graph $\Gamma$ of $G$. Define the \emph{coned off Cayley graph} of $G$ relative to $H_1,\ldots,H_m$ by adding a new vertex to $\Gamma$ for each left coset $gH_i$ of each $H_i$, and adding new edges connecting the new vertex to all the vertices of $gH_i$. We say that $\Gamma$ is \emph{weakly hyperbolic} relative to $H_1,\ldots,H_m$ if the coned off Cayley graph is a Gromov hyperbolic metric space. 

Choose unique representatives $[c_1],\ldots,[c_m]$ for the orbits of the action of $\MCG$ on $\C$, and let $H_i = \Stab[c_i] \subgroup \MCG$. It is easy to see \cite{MasurMinsky:complex1} that $\CC$ is equivariantly quasi-isometric to the coned off Cayley graph of $\MCG$ relative to $H_1,\ldots,H_m$. Much deeper is the question of what, exactly, is the large scale geometric behavior of these two spaces:

\begin{theorem}[\cite{MasurMinsky:complex1}]
The curve complex $\CC$ is a Gromov hyperbolic metric space, and $\MCG$ is weakly hyperbolic relative to stabilizers subgroups of representatives of orbits in $\C$.
\end{theorem}

The boundary of the curve complex was identified in the following theorem of Klarreich \cite{Klarreich:thesis}; see also \cite{Hamenstadt:ccboundary}. Two points in $\PMF$ are \emph{topologically equivalent} if they are represented by singular foliations which, after forgetting measure, are identical.

\begin{theorem}
\label{TheoremCCBoundary}
There is an $\MCG$ equivariant quotient map from the subspace of $\PMF$ consisting of classes of arational measured foliations to the Gromov boundary of $\CC$, such that the decomposition elements of this quotient map are the topological equivalence classes of arational measured foliations.
\end{theorem}

Using the fact that $\CC$ is convex cocompact, it makes sense to compare convex cocompactness of a subgroup of $\MCG$ with respect to its actions on $\Teich$ and on $\CC$. This is accomplished in the following theorem proved independently, and with different techniques, by Kent-Leininger and Hamenstadt:

\begin{theorem}[\cite{KentLeininger:Shadows}, \cite{Hamenstadt:extensions}]
\label{TheoremCurveComplex}
A finitely generated subgroup $G \subgroup \MCG$ is convex cocompact if and only if its action on $\CC$ is convex cocompact.
\end{theorem}

\paragraph{Schottky and virtual Schottky subgroups of mapping class groups.} 
\index{Schottky subgroup}
A convex cocompact, free subgroup of $\MCG$ is called a \emph{Schottky subgroup}. The ubiquity of Schottky subgroups was first demonstrated by Farb and myself in \cite{FarbMosher:quasiconvex}, using a result of \cite{Mosher:hypbyhyp} which is to be explained below. 

\begin{theorem}[\cite{FarbMosher:quasiconvex}]
For any independent set of pseudo-Anosov mapping classes $\phi_1,\ldots,\phi_K \in \MCG(S)$, there exists an integer $B \ge 1$ such that for all integers $\beta_1,\ldots,\beta_K \ge B$, the mapping classes $\phi_1^{\beta_1},\ldots,\phi_K^{\beta_K}$ freely generates a Schottky subgroup of $\MCG(S)$.
\end{theorem}

Kent and Leininger \cite{KentLeininger:Shadows} and independently Hamenstadt \cite{Hamenstadt:extensions} discovered a more direct and natural setting for this theorem, using a ping-pong argument in the curve complex. Here is their proof.

The action of $\phi_k$ on $\PMF$ has source sink dynamics, with source and sink each being arational, and for $k=1,\ldots,K$ the $2K$ fixed points are pairwise topologically inequivalent, because they are pairwise distinct and each is uniquely ergodic. By Theorem~\ref{TheoremCCBoundary}, it follows that the action of each $\phi_k$ on $\bdy\CC$ has source sink dynamics. and so $\phi_k$ is a hyperbolic isometry on the Gromov hyperbolic metric space $\CC$, and that their fixed point sets in $\bdy\CC$ are pairwise disjoint. To conclude that $\phi_1^{\beta_1},\ldots,\phi_K^{\beta_K}$ freely generate a convex cocompact subgroup of $\MCG$, we apply the following folk theorem:

\begin{theorem}[Hyperbolic ping-pong] 
\index{ping-pong!hyperbolic}
If $X$ is a Gromov hyperbolic geodesic metric space, and if $\phi_1,\ldots,\phi_K$ are hyperbolic isometries of $X$ with pairwise disjoint fixed points in $\bdy X$, then there exists an integer $B \ge 1$ such that for all integers $\beta_1,\ldots,\beta_K \ge B$, the isometries $\phi_1^{\beta_1},\ldots,\phi_K^{\beta_K}$ freely generate a convex cocompact group of isometries. 
\end{theorem}

It seems hard to track down an early proof of hyperbolic ping-pong, although it was certainly known to the very earliest practitioners of Gromov hyperbolicity. Here is a sketch of a proof, pretty much the same as the proof given in \cite{KentLeininger:Shadows}.

One may choose $\lambda,\epsilon$ quasigeodesic lines $\gamma_1,\ldots,\gamma_K$ in $X$ that are axes for $\phi_1,\ldots,\phi_K$, respectively (if $X$ were proper then one could choose $\lambda=1$, $\epsilon=0$, in other words, geodesics). Fix a base point $p \in X$, and let $D$ be an upper bound for the distance from $p$ to $\gamma_1,\ldots,\gamma_K$. Also, let $L$ be a lower bound for the length of a fundamental domain of each $\phi_1,\ldots,\phi_K$ along its  axis $\gamma_1,\ldots,\gamma_K$. Fix $B$ and $\beta_1,\ldots,\beta_K \ge B$, let $G$ be the group generated by $\phi_1^{\beta_1},\ldots,\phi_K^{\beta_K}$. 

We know by topological ping-pong that if $B$ is sufficiently large then $G$ is freely generated by $\phi_1^{\beta_1},\ldots,\phi_K^{\beta_K}$. It suffices to show that for each word $w = w_1 \cdots w_m$ in the generators, a geodesic from $p$ to $w(p)$ stays uniformly close to the orbit $G \cdot p$. We describe a path $\gamma_w$ from $p$ to $w(p)$ as follows. Denote $w_i = \phi_{k_i}^{\beta_{k_i}}$. Starting from $p$, jump a distance at most $D$ onto the axis $\gamma_{k_1}$, travel along $\gamma_1$ for $\beta_1$ fundamental domains, thereby travelling a distance of at least $BL$, then jump a distance at most $D$ to the point $w_1(p)$. Next, jump a distance at most $D$ onto the axis $w_1(\gamma_{k_2})$, travel along this axis for $\beta_2$ fundamental domains, thereby travelling a distance of at least $BL$, then jump a distance at most $D$ to the point $w_1 w_2(p)$. Continuing in this manner we get a path $\gamma_w$ from $p$ to $w(p)$, which stays in a $D+BL$ neighborhood of $G \cdot p$. The path $\gamma_w$ is a $(BL+2D)$ local $\lambda,(\epsilon+2D)$ quasigeodesic. If $B$ is sufficiently large, it follows that $\gamma_w$ is a $\lambda,\epsilon'$ quasigeodesic for $\epsilon'$ depending only on the previous constants and the hyperbolicity constant of $X$ \cite{Cannon:TheoryHyp}. Now we use the fact that $\gamma_w$ has Hausdorff distance at most $E$ from any geodesic $\ell$ with the same endpoints $p$, $w(p)$, where $E$ depends only on $\lambda$, $\epsilon'$, and the hyperbolicity constant for $X$. It follows that $\ell$ stays within distance $E+D+BL$ of $G \cdot p$.

\paragraph{Virtual Schottky subgroups.} A \emph{virtual Schottky subgroup} of $\MCG$ is a subgroup that contains a Schottky subgroup with finite index. Clearly convex cocompactness is preserved under passage to finite index subgroups, so every virtual Schottky subgroup is convex cocompact. The following recent result, to appear in the dissertation of Honglin Min, is another example of a ping-pong argument in the mapping class group:

\begin{theorem}[\cite{Min:HyperbolicGraphs}]
If $A,B$ are finite subgroups of $\MCG$, if $\phi \in \MCG(S)$ is pseudo-Anosov, and if the virtual normalizer of $\phi$ has trivial intersection with $A$ and $B$, then for all sufficiently large $n$ the subgroups $A$ and $\phi^n B \phi^{-n}$ freely generate their free product in $\MCG$, and this subgroup is virtually Schottky. 
\end{theorem}

\paragraph{Hyperbolic extensions of surface groups} 

When the surface $S$ is closed and oriented, there is a canonical short exact sequence
$$1 \to \pi_1(S,p) \to \MCG(S-p) \to \MCG(S) \to 1
$$
See \cite{Birman:Braids}. 

Given a subgroup $H \subset \MCG(S)$, let $\Gamma_H$ be its preimage in $\MCG(S-p)$, and so we have an extension
$$1 \to \pi_1(S,p) \to \Gamma_H \to H \to 1
$$
Thurston's geometrization theorem (see \cite{Otal:fibered}) shows that if $H$ is an infinite cyclic, pseudo-Anosov subgroup of $\MCG(S)$ then $\Gamma_H$ is a word hyperbolic group. 

The following theorem was the prequel to the theory of convex cocompact groups. It is proved by an application of the Bestvina--Feighn combination theorem.

\begin{theorem}[\cite{Mosher:hypbyhyp}] If $\phi_1,\ldots,\phi_K$ are independent pseudo-Anosov elements of $\MCG$ then sufficiently high powers of $\phi_1,\ldots,\phi_K$ freely generate a subgroup $H$ such that $\Gamma_H$ is word hyperbolic. 
\end{theorem}

With this result as motivation, in \cite{FarbMosher:quasiconvex} Farb and I set out to develop a geometric understanding of when $\Gamma_H$ is word hyperbolic, attempting to prove that for this to be true it was necessary and sufficient that $H$ be convex cocompact. We proved necessity, and when $H$ is free we proved sufficiency by using the Bestvina--Feighn combination theorem \cite{BestvinaFeighn:combination}. A completely general proof of sufficiency was recently given by Hamenstadt \cite{Hamenstadt:extensions}, so that now one can state:

\begin{theorem}[\cite{FarbMosher:quasiconvex}, \cite{Hamenstadt:extensions}]
\label{TheoremHypExtension}
For any subgroup $H \subset \MCG$, the group $\Gamma_H$ is word hyperbolic if and only if $H$ is convex cocompact.
\end{theorem}

\section{\Teichmuller\ discs and their stabilizers}
\label{SectionStabilizers}

One very well studied class of subgroups of $\MCG$ is stabilizers of \Teichmuller\ discs. A \emph{\Teichmuller\ disc}\index{Teichmuller disc} in $\Teich$ can be defined in one of the following equivalent ways. 

First, consider $\sigma \in \Teich$ and $\theta \in \QD_\sigma$. As $\omega$ varies over the unit circle in the complex plane, consider the family of quadratic differentials $\omega\theta$. This defines a family of rays based at $\sigma$ parameterized by the circle, and the union of these rays is defined to be the \emph{\Teichmuller\ disc with basepoint $\sigma$ and direction $\theta$}, denoted $D(\sigma,\theta)$. A somewhat more invariant way to say almost the same thing is that $D(\sigma,\theta)$ is the set of points in $\Teich$ obtained by choosing an element of $\PSL(2,\reals)$ acting affinely on the Euclidean plane, and applying that element to the regular canonical coordinates for $\theta$ to get a new conformal structure.

Second, a \Teichmuller\ disc $D$ is a maximal subset which is the image of a holomorphic embedding of the \Poincare\ disc into $\Teich$. A theorem of Royden \cite{Royden:AutIsomTeich} in the compact case and Gardiner \cite{Gardiner:approximation} in the noncompact case says that any two points $\sigma \ne \sigma' \in \Teich$ lie on a unique such disc $D$, and moreover that the \Teichmuller\ distance between $\sigma$ and $\sigma'$ equals their distance in the \Poincare\ metric on $D$. The identification of $D$ with $D(\sigma,\theta)$ for some $\theta$ is obtained by choosing $\theta$ to be the initial quadratic differential from $\sigma$ to~$\sigma'$. 

To each \Teichmuller\ disc $D$ there is associated an embedded circle $S^1(D) \subset \PMF$, defined as follows. Choose $\sigma \in \Teich$ and $\theta \in \QD_\sigma$ so that $D=D(\sigma,\theta)$. As $\omega$ varies over the unique circle in the complex plane, the horizontal measured foliations of $\omega\theta$ trace out the desired circle $S^1(D)$. Alternatively, a point $P[\fol] \in \PMF$ is in $S^1(D)$ if and only if it is represented as a linear measured foliation in any regular canonical coordinate $z$ for $\theta$; in other words, there exists an $\reals$-linear differential form $a \, dx + b \, dy$ such that $P[\fol]$ is represented in any $z$ by $\abs{a\, dx + b \, dy}$. This makes it clear that $S^1(D)$ is well-defined independent of $\sigma \in D$, because any other $\sigma' \in D$ is obtained by a \Teichmuller\ deformation with initial quadratic differential $\omega\theta$ for some $\omega$ in the unit circle, resulting in a terminal quadratic differential $\theta'$ on $\sigma'$, and linearity with respect to $\theta$ and to $\theta'$ are clearly equivalent.

Let $D \subset \Teich$ be a \Teichmuller\ disc, and let $\Stab(D) = \{\phi \in \MCG \suchthat \phi(D)=d\}$ be its stabilizer. The group $\Stab(D)$, or its image in $\PSL(2,\reals)$ when regarding $D$ as the upper half plane, is sometimes called a \emph{Veech group}, and the quotient orbifold $D / \Stab(D)$, a proper holomorphic curve in the moduli space of $S$, is sometimes called a \emph{Veech surface}. These terms are sometimes restricted to the case where $\Stab(D)$ is a lattice, meaning that $D / \Stab(D)$ has finite hyperbolic area.

Regarding $D$ as the hyperbolic plane, the elliptic--parabolic--loxodromic trichotomy for $\Isom(\hyp^2)$ becomes not just an analogy but a strict correspondence with the finite order--reducible--pseudo-Anosov dichotomy for elements of $\Stab(D)$; for details see see for example \cite{Veech:TeichmullerCurves}. To be precise, fix $\phi \in \Stab(D)$. Then $\phi$ is loxodromic on $D$ if and only if it is pseudo-Anosov in $\MCG$, $\phi$ is parabolic on $D$ if and only if it is reducible in $\MCG$, and $\phi$ is elliptic on $D$ if and only if it is finite order in $\MCG$. The case where $\phi$ is parabolic can be analyzed more closely. In this case, letting $D=D(\sigma,\theta)$, the action of $\phi$ is a shear transformation on any regular canonical coordinate $z$, parallel to the leaves of some constant slope measured foliation $\fol$ called the \emph{shear foliation} for $\phi$; this foliation corresponds to the point of $S^1(D)$ fixed by $\phi$. Each leaf of $\fol$ is compact, and so $\fol$ is obtained by enlarging some weighted family of simple closed curves $c_1,\ldots,c_n$. In particular, the essential curve family $\{c_1,\ldots,c_n\}$ is a canonical reducing system for $\phi$, and the first return of $\phi$ to each complementary component of $S-(c_1 \union \cdots \union c_n)$ is of finite order, so $\phi$ is the composition of a finite order mapping class that is reduced by $\{c_1,\ldots,c_n\}$ and a product of Dehn twists along the components of $c_1,\ldots,c_n$. 

Using this description, we can get a fairly precise understanding of the case when $\Stab(D)$ is a reducible subgroup of $\MCG$. Assuming this is so, let $\phi \in \Stab(D)$ be of infinite order, with shear foliation $\fol$. For any other infinite order $\phi' \in \Stab(D)$, the shear foliation of $\phi'$ must have the same slope as $\fol$ and so represent the same point in $S^1(D)$. Otherwise, if the shear foliation $\fol'$ of $\phi'$ is distinct from $\fol$ in $S^1(D)$, then high powers of $\phi$ and $\phi'$ would generate a loxodromic element of $\Stab(D)$, by the proof of the Tits alternative for $\hyp^2$ (this is not hyperbolic ping-pong, because $\phi$ and $\phi'$ do not have source--sink dynamics, but a variation of hyperbolic pong-pong works). Thus, every infinite order element of $\Stab(D)$ stabilizes the same point in $S^1(D)$, in fact every such element acts as a shear transformation on $\theta$ with respect to the same linear measured foliation $\fol_\omega$. But this implies that $\Stab(D)$ is virtually cyclic, generated by a multiple Dehn twist about components of the curves $c_1,\ldots,c_n$ whose enlargement is $\fol_\omega$. In this case we can identify the \McPapa\ limit set of $\Stab(D)$: it is the simplex in $\PMF$ spanned by the projective classes of $c_1,\ldots,c_n$.

Henceforth we shall only be interested in the case when $\Stab(D)$ is infinite and irreducible. In this case, since $S^1(D)$ is invariant under $\Stab(D)$ by construction, and since $S^1(D)$ is obviously closed in $\PMF$, by applying Theorem~\ref{TheoremMcPapaIrred} we obtain:

\begin{theorem}
\label{TheoremVeechLimits}
If $\Stab(D)$ is infinite and irreducible then the \McPapa\ limit set of $\Stab(D)$ is contained in $S^1(D)$.
\end{theorem}

Thurston \cite{Thurston:surfaces} gave the first examples of \Teichmuller\ discs $D$ such that $\Stab(D)$ is a lattice, meaning that it acts on $D$ with cofinite area. Start with the linear action of $\SL(2,\Z)$ on $T^2 = \reals^2 / \integers^2$, then choose a finite subset $A \subset T^2$ invariant under this action, and let $S \to T^2$ be a branched cover of $T^2$ branched over $A$. Lifting the conformal structure from $T^2$ to $S$ defines $\sigma \in \Teich(S)$, and lifting the quadratic differential $(dx + i\, dy)^2$ from $T^2$ to $S$ defines $\theta \in \QD_\sigma$. There exists a finite index subgroup $G \subgroup \SL(2,\Z)$ that lifts to $S$, and $G$ stabilizes $D(\sigma,\theta)$ and acts with cofinite area. Combining with Theorem~\ref{TheoremVeechLimits} it follows that the \McPapa\ limit set of $G$ is the circle $S^1(D)$.

Veech \cite{Veech:TeichmullerCurves} gave the first examples where $\Stab(D)$ is a lattice that does not arise from Thurston's construction of branching over the torus. Veech also developed some of the general theory of $\Stab(D)$, for instance establishing that $\Stab(D)$ never acts cocompactly on $D$. 
When $\Stab(D)$ is a lattice, it follows that there are finitely many maximal parabolic subgroups of $\Stab(D)$ up to conjugacy, and the Veech surface $D / \Stab(D)$ is a finite area holomorphic curve in the moduli space of $S$.

McMullen \cite{McMullen:complexity} proved that if $S$ is the closed surface of genus~2, if $D \subset \Teich(S)$ is a \Teichmuller\ disc, and if $\Stab(D)$ contains a hyperbolic element, then the limit set of $\Stab(D)$ is all of $S^1(D)$. In the same paper McMullen also gave examples where $\Stab(D)$ is infinitely generated, a phenomenon independently discovered by Hubert and Schmidt \cite{HubertSchmidt:veech}.

For a further survey of stabilizers of \Teichmuller\ discs, and references to the literature, see the section ``Veech surfaces'' in the problem list \cite{HMSZ:problems}, 
and the chapters by Harvey \cite{Harvey:dessins} and by Herrlich and Schmith\"usen \cite{HerrlichSchmithusen:Teichmuller} in this handbook.

\section{Surface subgroups and the Leininger--Reid combination theorem}
\label{SectionCombinations}
\index{combination theorems}
The only known examples of convex cocompact subgroups are free groups and virtually free groups. Even in $\Isom(\hyp^n)$, constructing nonfree discrete subgroups can be subtle. One of the pioneering tools in this regard is the Maskit combination theorem \cite{Maskit:Kleinian}.

One of the most interesting questions about subgroups of mapping class groups is whether there exists closed surfaces $S,F$ of genus $\ge 2$ and a convex cocompact subgroup of $\MCG(S)$ isomorphic to $\pi_1(F)$. If such a subgroup existed then, by Hamenstadt's theorem, one would obtain the first example of a compact \nb{4}manifold $M$ with word hyperbolic fundamental group that fibers over a surface. The preprint \cite{MKapovich:ComplexSurfaces} of M.\ Kapovich contains a theorem that says if such an $M$ existed then $M$ could not have a $\complex\hyp^2$ structure.

It is even hard to construct any surface subgroup $\pi_1(F)$ of $\MCG(S)$ with the genus of $F$ at least two. The first examples were given by Harvey and Gonz\'alez-D\'iez \cite{HGD:SurfaceGroups}, but they have many non-pseudo-Anosov elements and so cannot be convex cocompact.

The closest approach yet to a convex cocompact surface subgroup of $\MCG$ is the following theorem of Leininger and Reid, whose statement is reminiscent of the Maskit combination theorem. 

The starting point of the theorem is the construction of a \Teichmuller\ disc $D$ and finite index subgroup $G \subgroup \Stab(D)$ such that $\Stab(D)$ (and hence also $G$) is a lattice, $G$ is torsion free, and $G$ has exactly one parabolic subgroup, generated by a multitwist $\tau$ about a multicurve $c$ such that no component of $S-c$ is a three holed sphere. It follows that there exists a homeomorphism $\psi$ that preserves $c$ and is pseudo-Anosov on each component of $S-c$.

\begin{theorem}[\cite{LeiningerReid:combination}]
With the objects as denoted above, for each sufficiently large $k$ the group generated by $G$ and $\phi^k G \phi^{-k}$ is a free product with amalgamation along $\<\tau\>$, and is therefore isomorphic to $\pi_1(F)$ for an oriented surface $F$ of genus $\ge 2$. Moreover, every element of $\pi_1(F) \subgroup \MCG(S)$ is pseudo-Anosov \emph{except} for elements conjugate to powers of $\tau$.
\end{theorem}

A Leininger--Reid subgroup cannot be convex cocompact, because it contains a multitwist. In \cite{Mosher:ExtensionProblems}, a notion of geometric finiteness for subgroups of $\MCG$ is proposed, and the question is posed whether Leininger--Reid subgroups are geometrically finite.

%%%%%%%%%%%%%%%%%%%%%%%%%%%%%%%%%
%\bibliographystyle{amsalpha} %
%\bibliography{mosher}        %
%%%%%%%%%%%%%%%%%%%%%%%%%%%%%%%%%

\newcommand{\etalchar}[1]{$^{#1}$}
\def\cprime{$'$} \def\cprime{$'$}
\providecommand{\bysame}{\leavevmode\hbox to3em{\hrulefill}\thinspace}
\providecommand{\MR}{\relax\ifhmode\unskip\space\fi MR }
\providecommand{\MRhref}[2]{%
  \href{http://www.ams.org/mathscinet-getitem?mr=#1}{#2}
}
\providecommand{\href}[2]{#2}

\end{document}